\documentclass[10pt]{amsart}   
\usepackage{amssymb}
\usepackage{amsthm}
\usepackage{amsmath}
\usepackage{amscd}
\input xy
\xyoption{all}

\newtheorem{thm}{Theorem}[section]
\newtheorem{prop}[thm]{Proposition}
\newtheorem{clly}[thm]{Corollary}
\newtheorem{lemma}[thm]{Lemma}

\theoremstyle{definition}
\newtheorem{defi}[thm]{Definition}
\theoremstyle{example}
\newtheorem{ex}[thm]{Example}

\newtheorem{nota}[thm]{Notation}

\theoremstyle{Remark}
\newtheorem{rk}[thm]{Remark}

\newcommand{\KK}{{\mathcal{K}}}
\newcommand{\LL}{{\mathcal{L}}}
\newcommand{\NN}{{\mathbb{N}}}
\newcommand{\ZZ}{{\mathbb{Z}}}

\newcommand{\CC}{{\mathbb{C}}}

\newcommand{\cstar}{\mbox{$C^*$}}

\newcommand{\id}{\mbox{id}}

\newcommand{\ox}{\mbox{${\tilde{{\mathcal O}}}_X$}}
\newcommand{\xna}{A_n^X}
\newcommand{\xn}{X_n}
\newcommand{\xnm}{M_n^X}

\newcommand{\ainf}{A_{\infty}^X}
\newcommand{\xinf}{X_\infty}

\newcommand{\rxa}{X(A)}
\newcommand{\rxm}{X(M)}
\newcommand{\rx}{X(X)}

\begin{document}

\title[Cuntz-Pimsner \cstar-algebras and crossed products]{Cuntz-Pimsner \cstar-algebras and crossed products by Hilbert \cstar-bimodules}
\author{Beatriz Abadie}
\thanks{Partially supported by Proyecto Conicyt Clemente Estable 8013}

\address{Centro de Matem\'aticas, Facultad de Ciencias, 
Igu\'a 4225, CP 11 400, Montevideo, Uruguay.}
\email{abadie@cmat.edu.uy}
\author{Mauricio Achigar}
\address{Centro de Matem\'aticas, Facultad de Ciencias, 
Igu\'a 4225, CP 11 400, Montevideo, Uruguay.}
\email{achigar@cmat.edu.uy}

\date{\today}

\begin{flushright}
\end{flushright}
\vspace{1cm}

\begin{abstract}
Given a correspondence $X$ over a \cstar-algebra $A$, we construct a \cstar-algebra $A^X_{\infty}$ and a Hilbert \cstar-bimodule $X_{\infty}$ over   $A^X_{\infty}$ such that the augmented Cuntz-Pimsner \cstar-algebras ${\mathcal{\tilde{ O}}}_X$ and the crossed product $A^X_{\infty}\rtimes X_{\infty}$ are isomorphic. This construction enables us to  establish a condition for two augmented Cuntz-Pimsner \cstar-algebras to be Morita equivalent.
\end{abstract}

\maketitle

\section{Introduction and preliminaries}

The augmented Cuntz-Pimsner \cstar-algebra $\ox$ defined in \cite{pim} is a \cstar-algebra associated to an $A$-correspondence $(X,\phi_X)$  that is universal for certain  covariance conditions (see \cite[3.12]{pim}) when $\phi_X$ is injective and $X$ is full as a right Hilbert \cstar-module. 

On the other hand, when $X$ is also a Hilbert \cstar-bimodule over the \cstar-algebra  $A$  the crossed product $A\rtimes X$ defined in \cite {aee} is universal for  covariance conditions that  agree with those for which $\ox$ is universal under the assumptions mentioned above. 

Thus both constructions can  be carried out when $X$ is a Hilbert \cstar-bimodule, and they agree when $X$ is full on the right and the action on the left is faithful. But this may fail if the condition of faithfulness of the left action is dropped,  as the following example, shown to us by S{{\o}}ren Eilers, proves.

Let  $A=\CC\oplus \CC$ and $X=\CC$ be the Hilbert \cstar-bimodule over $A$ obtained by setting: 
\[(\lambda,\mu)\cdot x=\lambda x,\ x\cdot (\lambda,\mu)=x\mu,\  \langle x,y\rangle_L=(x\overline{y},0)\mbox{  and  }\langle x,y\rangle_R=(0,\overline{x}y).\]
Then $X\otimes X=0$ because $x\otimes y=x\cdot(0,1)\otimes y=x\otimes (0,1)\cdot y=x\otimes 0$. This implies that $\ox=\{0\}$ whereas  $A\rtimes X$ is isomorphic to $M_2(\CC)$. This last statement  can be checked directly by verifying that the *-homomorphism  induced by the covariant pair of maps $(i_A,i_X):(A,X)\rightarrow M_2(\CC)$ given by $i_A(\lambda,\mu)=\left(\begin{smallmatrix}
\lambda & 0\\
0 & \mu
\end{smallmatrix}\right)$, $i_X(x)=\left(\begin{smallmatrix}
0 & x\\
0 & 0
\end{smallmatrix}\right)$ is an isomorphism, or by noting that $X$ is the bimodule associated (as described in \cite[3.2]{aee}) to the partial action on $A$ given by $I=\CC\oplus 0$, $J=0\oplus \CC$, $ \theta(x,0)=(0,x)$.

This shows that $A\rtimes X$ and $\ox$ may not agree for a Hilbert \cstar-bimodule $X$ over $A$. 
On the other hand, as mentioned in  \cite[1]{aee}, for any $A$-correspondence $(X,\phi_X)$ the algebra $\ox$ is  a crossed product $A_\infty\rtimes X_\infty$. The example above shows that  the algebra $A_{\infty}$ and the bimodule $X_{\infty}$ do  not necessarily agree with the original $A$ and  $X$  when $X$ is a  Hilbert \cstar- bimodule over $A$. In this work we give an abstract construction of $A_\infty$ and $X_{\infty}$,  out of  an $A$-correspondence $(X,\phi_X)$. Both  $A_\infty$ and $ X_{\infty}$ are described as direct limits of nicely related directed sequences in their respective categories.

We apply this construction to the discussion of Morita equivalence of augmented Cuntz-Pimsner \cstar-algebras. One of our tools is a  result from  \cite[4.2]{aee}: if $X$ and $Y$ are Hilbert \cstar-bimodules over \cstar-algebras $A$ and $B$  respectively, and $M$ is an $A-B$ Morita equivalence bimodule such that the $A-B$ Hilbert \cstar bimodules $X\otimes_AM$ and $M\otimes_BY$ are isomorphic, then the crossed products $A\rtimes X$ and $B\rtimes Y$ are Morita equivalent.
In Theorem \ref{meq} we establish a condition of this kind for the Morita equivalence of two  augmented Cuntz-Pimsner  \cstar-algebras.   
P. Muhly and B. Solel showed in \cite[3.3, 3.5]{ms} a similar   result for  Cuntz-Pimsner \cstar-algebras and for correspondences $(X,\phi_X)$ and $(Y,\phi_Y)$ such that the maps $\phi_X$ and $\phi_Y$ are injective and the correspondences are non-degenerate (that is, $\phi_X(A)X=X$ and similarly for $Y$). 
 Our result for  augmented \cstar-algebras   does not require the action to be injective, but  a condition related to non-degeneracy (see  Remark \ref{muso}) has to be met.

This work is organized as follows. Section 2 deals with the notion of direct limit of Hilbert \cstar-modules and proves some basic results that will be further required. In section 3 we construct, for an $A$-correspondence $(X,\phi_X)$,  a \cstar-algebra $A_\infty$ and a Hilbert \cstar-bimodule $X_\infty$ over $A_\infty$ such that $\ox$ and $A_\infty \rtimes X_\infty$ are isomorphic. In section 4 we use that construction together with  \cite[4.2]{aee} to give a sufficient condition for the Morita equivalence of two augmented Cuntz-Pimsner \cstar-algebras.

We  start by recalling some definitions and by setting some notation.

\begin{nota}
\label{not}
{\rm Let $A$ and $B$ be \cstar-algebras. If  $\phi:A\longrightarrow B$ is   a *-homomorphism, we denote by $\phi^{(k)}$ the *-homomorphism  $\phi^{(k)}:M_k(A)\longrightarrow M_k(B)$ defined by  $(\phi^{(k)}(M))_{ij}=\phi(M_{ij})$.

Given a right Hilbert \cstar-module $X$ over a \cstar-algebra $A$ we  denote by $\LL(X)$ and $\KK(X)$, respectively, the \cstar-algebras of adjointable and compact maps. For $x,y\in X$, we write $\theta_{x,y}$ to denote the map  $\theta_{x,y}\in \KK(X)$ defined by  $\theta_{x,y}(z)=x\langle y, z\rangle$. For $x\in X$,  $|x|$ denotes the element  $|x|\in A$ defined by $|x|=\langle x,x\rangle ^{1/2}$.

Given subsets $S$ and $T$ of $X$, we write $\langle S,T\rangle$ to denote the set $\langle S,T\rangle =\overline{\rm{span}}\{\langle s, t\rangle :s\in S,\ t\in T\}$.
If $S\subset \LL(X)$,  we denote by $SX$ the set $SX=\overline{\rm{span}}\{s(x):s\in S,\ x\in X\}$. 
Given a \cstar-subalgebra $C$ of $\LL(X)$,  we denote by $L_{C,X}$  the right Hilbert \cstar-module homomorphism $L_{C,X}:C\otimes_C X\longrightarrow X$ defined  by $L_{C,X}(c\otimes x)=c(x)$. Note that  $L_{C,X}$  is an isomorphism when $CX=X$.

When $X$ is a right Hilbert \cstar-module over $A$, the map $x\otimes a\mapsto xa$, for $x\in X$ and $a\in A$,  is an isomorphism of right A-Hilbert \cstar-modules between $X\otimes A$ and $X$ that associates the map $T\in \LL(X)$ to the map $T\otimes \id_A\in\LL(X\otimes A)$. Often in  this work we will  identify  $X$ with $X\otimes A$ and $T\in \LL(X)$ with $T\otimes \id_A$ as above without further warning. 
 For $T\in \LL(X)$ we will understand  that $T^{\otimes 0}$ is $\id_A$.

We next recall some of the terminology in \cite{ms} that we will adopt. Given \cstar-algebras $A$ and $B$, an $A-B$ correspondence $(X,\phi_X)$ consists of a right Hilbert \cstar-module $X$ over $B$ together with  a $*$-homomorphism $\phi_X:A\longrightarrow \LL(X)$. We will  denote the correspondence by $X$ and drop the reference to the map $\phi_X$ when it does not lead to confusion. Besides, we will   write $a\cdot x$  to denote $[\phi_X(a)](x)$.

Let $X_i$ be an $A_i-B_i$ correspondence, for $i=1,2$. A homomomorphism  of correspondences $(\sigma,\phi,\pi)$ consists of \cstar-algebra homomorphisms $\sigma:A_1\longrightarrow A_2$ and $\pi:B_1\longrightarrow B_2$ and a linear map $\phi:X_1\longrightarrow X_2$ such that 
\[\phi (a\cdot xb)=\sigma(a)\cdot \phi(x)\pi(b)\  \mbox{ and }\ \langle \phi(x),\phi(y)\rangle=\pi\big(\langle x,y\rangle\big),\]
for all $x,y\in X_1,\ a\in A_1,\ b\in B_1.$ 

Whenever $A_1=A_2$ (respectively $B_1=B_2$) and there is no reference to the map  $\sigma$ (respectively $\pi$) we assume it is the identity map. 
 Two $A-B$ correspondences $X$ and $Y$ are said to be isomorphic if there is a homomorphism $(\id_A,J,\id_B)$ where  $J:X\longrightarrow Y$ is invertible.

Note that the map $L_{C,X}$ defined above for a \cstar-subalgebra $C$ of $\LL(X)$ is a homomorphism of $C-A$ correspondences, for $C$ acting on $C\otimes X$ via  $l\otimes\id_X$, $l$ being left multiplication.

Homomorphisms of right Hilbert \cstar-modules and homomorphisms of  Hilbert \cstar-bimodules are defined in the obvious analogous way, Hilbert \cstar-bimodules being defined as in \cite[1.8]{bms}.}
\end{nota}
 
\begin{lemma}
\label{norm}
Let $(\phi, \sigma): (X,A)\longrightarrow (Y,B)$ be a homomorphism of right Hilbert \cstar-modules. Then $\phi$ is norm-decreasing,  and it induces a \cstar-algebra homomorphism $\phi_*:\KK(X)\longrightarrow \KK(Y)$ such that $\phi_*(\theta_{x_1,x_2})=\theta_{\phi(x_1),\phi(x_2)}$, for  $x_1$,  $x_2\in X$.

\end{lemma}
\begin{proof}

If $x\in X$, then
\[\|\phi(x)\|^2=\|\langle \phi(x),\phi(x)\rangle\|=\|\sigma(\langle x,x\rangle)\|\leq \|\langle x,x\rangle\|=\|x\|^2.\]
As for the second statement, if $x_i,y_i\in X$ for $i=1,2,\dots,n$, then by \cite[2.1]{pin} we have 
\[\|\sum \theta_{\phi(x_i),\phi(y_i)}\|=\|S^{1/2}T^{1/2}\|_{M_n(B)},\] 
where $S_{ij}=\langle \phi(x_i),\phi(x_j)\rangle$ and $T_{ij}=\langle \phi(y_i),\phi(y_j)\rangle$.
Now,  $S= \sigma^{(n)}(M)$ and $T= \sigma^{(n)}(N)$, where $M_{ij}=\langle x_i,x_j\rangle$ and $N_{ij}=\langle y_i,y_j\rangle$. Therefore 
\[\|S^{1/2}T^{1/2}\|=\|\sigma^{(n)}\big(M^{1/2}N^{1/2}\big)\|\leq \|M^{1/2}N^{1/2}\|=\|\sum \theta_{x_i,y_i}\|,\]
which shows that $\phi_*$ extends to a continuous map on $\KK(X)$. Finally,  it is straightforward to check that $\phi_*$ is a $*$-homomorphism from the fact that $(\sigma,\phi)$ is a homomorphism of right Hilbert \cstar-modules.  

\end{proof}

\section{Directed sequences of right Hilbert \cstar-modules}

In this section we discuss a procedure to get, for a given $A$-correspondence $(X,\phi_X)$,  a Hilbert \cstar-bimodule $X_\infty$ over a \cstar-algebra $A_{\infty}$. We will show in next section that  $A_\infty\rtimes  X_\infty$ is isomorphic to $\ox$.

In order to get a left inner product on $X$ one needs to add to Im $\phi$ the compact operators $\KK(X)$. If one lets $A_1\subset \LL(X)$ be defined by $A_1=$ Im $\phi +\KK(X)$, then $X$ is an $A_1-A$ Hilbert \cstar-bimodule, but there is no clear right action of $A_1$ on $X$. This suggests replacing $X$ by $X_1:=X\otimes_A A_1$. Thus we end up  with an $A_1$-correspondence $X_1$,  and the procedure can be iterated. We show how this iteration yields directed sequences $\{A_n\}$ and $\{X_n\}$ whose limits  $A_\infty$ and $X_{\infty}$ are such that  $X_{\infty}$ is a Hilbert \cstar-bimodule over $A_\infty$. We will develop this procedure in a somewhat more general context that will be of use in the discussion of Morita equivalence in the last section.

\begin{defi}

A directed sequence  $\{(X_n,A_n,\phi_n^X, \phi_n^A)\}$ of right Hilbert \cstar-modules consists of a directed sequence $\{(A_n,\phi_n^A)\}$ of \cstar-algebras together with   a directed sequence $\{(X_n,\phi_n^X)\}$ of vector spaces such that $X_n$ is a right Hilbert \cstar-module over $A_n$ and $(\phi_n^X, \phi_n^A)$ is a homomorphism of right Hilbert \cstar-modules for each $n\in \NN$.

\end{defi}

\begin{rk}

{\rm{Let $\{(X_n,A_n,\phi_n^X, \phi_n^A)\}$ be a  directed sequence   of right Hilbert \cstar-modules.  Since  the maps $\phi_n^X$ are norm decreasing by Lemma \ref{norm}, the sequence  $\{X_n,\phi_n^X\}$ has a direct  limit $(X_\infty, \{\lambda_n^X\})$,  that can be described as follows. 
Let $Y_0$ be the vector space 
\[Y_0=\{x\in\prod X_n: \mbox{ there exists }n_0 \mbox{ such that } x_{n+1}= \phi_n(x_n) \  \forall n\geq n_0\},\]
 and let $Y=\{x\in Y_0: \lim_n\|x_n\|=0\}$. Then $X_\infty$ is the completion of $Y_0/Y$ for the norm $\|x\|=\lim_n\|x_n\|$. The canonical maps $\lambda_n^X:X_n\longrightarrow X_\infty$ are given by $\lambda_n^X=\pi\circ\tilde{\lambda}_n^X$, where $\pi:Y_0\longrightarrow X_\infty$ is the canonical projection and $\tilde{\lambda}_n^X(x_n)(k)=\phi^X_{n,k}(x_n)$, for $\phi^X_{n,k}: X_n\to X_k$  given by:

\[\phi^{X}_{n,k}=
\begin{cases}
0 &\text{if $k<n$;}\\
\id&\text{if $k=n$;}\\
\phi^{X}_{k-1}\circ\phi^{X}_{k-2}\dots \circ\phi^{X}_n &\text{if $k>n$.}
\end{cases}\]
Note that $\|\lambda_n^X(x)\|_{X_\infty}=\lim_m \|\phi^X_{n,m}(x_n)\|_{X_m}=\inf_m \|\phi^X_{n,m}(x_n)\|_{X_m}.$

If $\{x_n\}\in Y_0$,  and $n_0$ is such that  $x_{n+1}= \phi_n(x_n)$ for all $n\geq n_0$, then $\lambda_{n_0}^X(x_{n_0})=\pi\big(\{x_n\}\big)$, which shows that $\bigcup \lambda^X_n(X_n)$ is dense in $X_\infty$.

It is well known that  a similar description holds for the direct limit  $(A_{\infty}, \{\lambda_n^A\})$ of the directed sequence   of \cstar-algebras $\{A_n,\phi_n^A\}$, and that $\bigcup \lambda^A_n(A_n)$ is dense in $A_\infty$.  Note that $(\phi^A_{n,k},\phi^X_{n,k})$ is a homomorphism of right Hilbert \cstar-modules  from $(A_n,X_n)$ to $(A_k,X_k)$.

We will say that  $(X_\infty,A_\infty,\{ \lambda_n^X\}, \{ \lambda_n^A\} )$ is the direct limit of the directed sequence  $\{(X_n,A_n,\phi_n^X, \phi_n^A)\}$.}}

\end{rk}

\begin{prop}
\label{dlim}
 Let   $\{(X_n,A_n,\phi_n^X, \phi_n^A)\}$  be  a directed sequence   of right Hilbert \cstar-modules with direct limit $(X_\infty,A_\infty,\{\lambda_n^X\}, \{\lambda_n^A\})$. Then $X_\infty$ can be made into a right Hilbert \cstar-module over $A_{\infty}$ by setting:
\[\lambda^X_n(x_n)\lambda^A_n(a_n):=\lambda^X_n(x_na_n),\ \langle \lambda^X_n(x^n_1), \lambda^X_n(x^n_2) \rangle:= \lambda^A_n(\langle x^n_1,x^n_2\rangle),\]
for  $a_n\in A_n$,  $x_n, x^n_i\in X_n$, $i=1,2$. (Therefore  $(\lambda^X_n,\lambda^A_n): (X_n,A_n)\longrightarrow (X_\infty,A_\infty)$ is a homomorphism of right Hilbert \cstar-modules for all $n$.) 

Let $M$ be a right Hilbert \cstar-module over a  \cstar-algebra $B$ and, for each $n\in \NN$, let $(\mu_n^X,\mu_n^A):(X_n, A_n)\longrightarrow (M,B)$ be a homomorphism of right Hilbert \cstar-modules, such that the diagrams
\[\begin{array}{rr}
\xymatrix{
X_n\ar[r]^{\phi_n^X} \ar[d]_{\mu^X_n}&X_{n+1}\ar[dl]^{\mu_{n+1}^X}\\
M}&\xymatrix{
A_n\ar[r]^{\phi_n^A} \ar[d]_{\mu^A_n}&A_{n+1}\ar[dl]^{\mu_{n+1}^A}\\
B}
\end{array}\]
commute. If $\mu^X:X_\infty\longrightarrow M$, $\mu^A:A_\infty\longrightarrow B$ are the canonical maps yielded by the universal property of the direct limit, then  $(\mu^X,\mu^A)$ is a homomorphism of Hilbert \cstar-modules.

Besides,  the norm on $X_\infty$ induced by its structure of  right  $A_\infty$-Hilbert \cstar-module  agrees with the original norm.
\end{prop}

\begin{proof} We first check that the definition of the action on the right  makes sense.  Assume that   $\lambda^A_k(a_k)=\lambda^A_n(a_n)$ and that  $\lambda^X_k(x_k)=\lambda^X_n(x_n)$,  for some $a_n\in A_n$,  $a_k\in A_k$, $x_n\in X_n$, and $x_k\in X_k$. 

Given   $\epsilon>0$,  choose $j\in \NN$, $j\geq k$, $j\geq n$, and  large enough to have $\|\phi^X_{n,j}(x_n)-\phi^X_{k,j}(x_k)\|<\epsilon$ and $\|\phi^A_{n,j}(a_n)-\phi^A_{k,j}(a_k)\|<\epsilon$. Then
\[\|\lambda^X_n(x_na_n)-\lambda^X_k(x_ka_k)\|\leq \|\phi^X_{n,j}(x_na_n)-\phi^X_{k,j}(x_ka_k)\|\]
\[=\|\phi^X_{n,j}(x_n)\phi^A_{n,j}(a_n)-\phi^X_{k,j}(x_k)\phi^A_{k,j}(a_k)\|\]
\[\leq \|\phi^X_{n,j}(x_n)(\phi^A_{n,j}(a_n)-\phi^A_{k,j}(a_k))\| + \|(\phi^X_{n,j}(x_n)-\phi^X_{k,j}(x_k))\phi^A_{k,j}(a_k)\|\]
\[\leq (\|x_n\|+\|a_k\|)\epsilon.\]
 Besides,
\[\|\lambda_n(x_na_n)\|=\lim_m\|\phi^X_{n,m}(x_na_n)\|\]
\[=\lim_m\|\phi^X_{n,m}(x_n)\phi^A_{n,m}(a_n)\|\]
\[\leq \big(\lim_m\|\phi^X_{n,m}(x_n)\|\big)\big(\lim_m\|\phi^A_{n,m}(a_n)\|\big)\]
\[=\|x\|\|a\|,\]
which shows that the right action of $\bigcup_n\lambda^A_n(A_n)$ on $\bigcup_n\lambda^X_n(X_n)$ extends by continuity to a right  action of $A_\infty$ on $X_{\infty}$.

As for the definition of the right inner product, it  makes sense because if   $\lambda^X_n(x^n_i)=\lambda^X_k(x^k_i)$ for some $x^n_i\in X_n$, $x^k_i\in X_k$, $i=1,2$, then for any  $\epsilon>0$  we can choose $j\in \NN$ such that $j\geq k$, $j\geq n$, and  $\|\phi^X_{n,j}(x^n_i)-\phi^X_{k,j}(x^k_i)\|<\epsilon$ for $i=1,2$.
Then:
\[\|\lambda^A_n(\langle x^n_1,x^n_2\rangle)-\lambda^A_k(\langle x^k_1,x^k_2\rangle)\|\leq\|\phi^A_{n,j}(\langle x^n_1,x^n_2\rangle)-\phi^A_{k,j}(\langle x^k_1,x^k_2\rangle)\|\]
\[=\|\langle \phi^X_{n,j}(x^n_1), \phi^X_{n,j}(x^n_2)\rangle-\langle \phi^X_{k,j}(x^k_1), \phi^X_{k,j}(x^k_2)\rangle\|\]
\[\leq \|\langle \phi^X_{n,j}(x^n_1)-\phi^X_{k,j}(x^k_1), \phi^X_{n,j}(x^n_2)\rangle\|+ \|\langle \phi^X_{k,j}(x^k_1), \phi^X_{n,j}(x^n_2) -  \phi^X_{k,j}(x^k_2)\rangle\|\]
\[<\epsilon \big(\|x^n_2\|+\|x^k_1\|\big).\] 
Also note that, for $x_n\in X_n$ we have 
\[\begin{array}{ll}
\|\langle \lambda_n^X(x_n),\lambda^X_n(x_n)\rangle \|_{A_\infty}&=\|\lambda_n^A\big(\langle x_n,x_n\rangle\big)\|_{A_\infty}\\
& =\lim_m\|\phi^A_{n,m}\big(\langle x_n,x_n\rangle\big)\|_{A_m}\\
&=\lim_m\|\langle \phi^X_{n,m}(x_n),\phi^X_{n,m}(x_n)\rangle\|_{A_m}\\
&=\lim_m\|\phi^X_{n,m}(x_n)\|^2_{X_m}\\
&=\|\lambda^X_n(x_n)\|^2_{X_\infty},
\end{array}\]
which shows that the two norms on $X_\infty$ agree.

The remaining properties and statements are apparent from the definitions.

\end{proof}
\begin{ex}
\label{main}
{\rm{The following example will be of importance in this work. Given a correspondence $(X,\phi_X)$ over a \cstar-algebra $A$,  let $\rxa$ denote the \cstar-subalgebra  $X(A)={\mathcal K}(X)+ \mbox{Im } \phi_X$.
 Note that $X$ is an $X(A)-A$ Hilbert \cstar-bimodule.

Given  a right $A$-Hilbert \cstar-module $M$, we define the  right $X(A)$-Hilbert \cstar-module $\rxm$ by $\rxm:=M\otimes_{\phi_X}\rxa$, where $\rxa$ is viewed as an $ \mbox{Im } \phi_X-\rxa$ correspondence in the obvious way.

Note that

\[\begin{array}{ll}
\langle \sum_{i=1}^n m_i\otimes \phi_X(a_i),\sum_{j=1}^m p_j\otimes \phi_X(b_j)\rangle  &=  \sum_{i,j}\langle m_i\otimes \phi_X(a_i), p_j\otimes \phi_X(b_j)\rangle\\
&=\sum_{i,j}\phi_X(a_i)^*\phi_X(\langle m_i, p_j\rangle)\phi_X(b_j)\\
&= \sum_{i,j}\phi_X(\langle m_ia_i,p_jb_j\rangle\\
&=\phi_X (\langle\sum_i m_ia_i,\sum_j p_jb_j\rangle),
\end{array}\]
for $m_i, p_j\in M$,  $a_i,b_j\in A$, and $i=1,2,\dots,n$, $j=1,2,\dots,m$  

In particular $\|\sum_{i=1}^n m_i\otimes \phi_X(a_i)\|^2=\|\phi_X (|\sum_i m_ia_i|)\|^2 \leq \|\sum_i m_ia_i\|^2.$
This shows that one can define a map  $\psi_M^X: M\longrightarrow \rxm$ by $\psi^X_M(ma)=m\otimes \phi_X(a)$ so that $(\psi_M^X,\phi_X)$ is a homomorphism of right Hilbert \cstar-modules.  When $M=X$ the map $\psi^X_X$  will  be denoted  by $\psi_X$. In this case  $(\psi_X,\phi_X)$ is a homomorphism of correspondences.

Now, since $\rx$ and  $X(M)$ are, respectively, a correspondence and  a right Hilbert \cstar-module over $\rxa$,  the construction above can be iterated to get a sequence $\{A^X_n\}_{n\geq 0}$ of \cstar-algebras and, for each $n\geq 0$,  a correspondence $X_n$ over $A^X_n$ and a right $A^X_n$-Hilbert \cstar-module $M^X_n$ by setting $A^X_0=A$,  $X_0=X$, $M^X_0=M$, and, for $n\geq 0$:
\[A_{n+1}^X=\xn(\xna),\ X_{n+1}=X_n(X_n),\ \mbox{and } M_{n+1}^X=X_n(\xnm).\]

We  also get right Hilbert \cstar-module homomorphisms
\[(\phi_n^{M,X},\phi_n^{A,X}):(M_n^X,\xna)\longrightarrow (M_{n+1}^X,A_{n+1}^X) \mbox{ for all $n\geq 0$,  }\]
given by $\phi_n^{A,X}=\phi_{X_n}$ and  $\phi_n^{M,X}=\psi^{X_n}_{M_n}$,  that is, $\phi_n^{A,X}(a)=a\otimes {\rm{id}}_{A_n^X}$ for all $n\geq 1$, and  $\phi_n^{M,X}(ma)=m\otimes \phi_n^{A,X}(a)$, for $m\in M_n^X$  and $a\in A_n^X$.

When $M=X$ we write $\phi_n^X$ in place of  $\phi_n^{M,X}$. In that case $( \phi_{n}^{A,X}, \phi_{n}^X, \phi_n^{A,X})$ and $( \phi_{n+1}^{A,X}, \phi_{n}^X, \phi_n^{A,X})$ are, respectively,   homomorphisms of   correspondences and  Hilbert \cstar-bimodules:
\[\begin{array}{ll}
\big(\theta_{\phi^X_n(xa), \phi^X_n(yb)}\big)(z\otimes c)&=x\otimes \phi_n^{A,X}(a)\langle y\otimes \phi_n^{A,X}(b), z\otimes c\rangle\\
&=x\otimes \phi_n^{A,X}(a)\phi_n^{A,X}(b^*\langle y,z\rangle)c\\
&=xa\langle yb,z\rangle\otimes c\\
&=\big(\phi_{n+1}^{A,X}(\theta_{xa,xb})\big)(z\otimes c),
\end{array}\]
for all $x,y,z\in X_n$, $a,b\in A_n$ and $c\in A_{n+1}$.

Let  $(X_\infty,A^X_\infty,\{\lambda_n^X\},\{\lambda_n^A\})$ and 
 $(M^X_\infty,A^X_\infty,\{\lambda_n^M\},\{\lambda_n^A\})$ denote the direct limits of the sequences  $\{(X_n,A^X_n,\phi_n^X,\phi_n^{A,X})\}$ and 
 $\{(M^X_n,A^X_n,\phi_n^{M,X},\phi_n^{A,X})\}$, respectively. By Proposition \ref{dlim},  both $X_\infty$ and $M^X_\infty$ are right Hilbert \cstar-modules over $A^X_\infty$. }}

\end{ex}
\begin{rk}
\label{bimodule}
{\rm{In fact,  $X_\infty$ is a Hilbert \cstar-bimodule over $\ainf$: since $X_n$ is an $A^X_{n+1}-A^X_{n}$ Hilbert \cstar-bimodule  for all $n\in\NN$, the proof of Proposition \ref{dlim} carries over to the left structure of $X_\infty$, and the compatibility between the left and the right structures on $X_\infty$ is easily checked.}}
\end{rk}

\begin{prop}
\label{comp}
 Let   $\{(X_n,A_n,\phi_n^X, \phi_n^A)\}$  be  a directed sequence   of right Hilbert \cstar-modules with direct limit $(X_\infty,A_\infty,\{\lambda_n^X\}, \{\lambda_n^A\})$. 

Then $\big(\KK(X_\infty), \{(\lambda_n^X)_*\}\big)$ is the direct limit of $\{\big(\KK(X_n),(\phi_n^X)_*\big)\}$, where  $(\lambda_n^X)_*$ and $(\phi_n^X)_*$ are  defined  as in Lemma \ref{norm}.

\end{prop}

\begin{proof}

It is well known that for any integer $k$,  $(M_k(A_{\infty}),\{ (\lambda^A_n)^{(k)}\})$ is the direct limit of $\{M_k(A_n), (\phi^A_n)^{(k)}\}$, which in particular implies 
 that 
\[\|(\lambda^A_n)^{(k)}(T)\|_{M_k(A_\infty)}=\lim_m\|\big(\phi_{n,m}^A\big)^{(k)}(T)\|\mbox{  for all }T\in M_k(A_n).\]

 Now, the commuting diagram
\[
\xymatrix{
X_n\ar[r]^{\phi_n^X} \ar[d]_{\lambda^X_n}&X_{n+1}\ar[dl]^{\lambda_{n+1}^X}\\
X_\infty}\]
yields a commuting diagram
\[
\xymatrix{
\KK(X_n)\ar[r]^{(\phi_n^X)_*} \ar[d]_{(\lambda^X_n)_*}&\KK(X_{n+1})\ar[dl]^{(\lambda_{n+1}^X)_*}\\
\KK(X_\infty)}\]
which in turn yields a map $H:\varinjlim\KK(X_n)\longrightarrow\KK(X_\infty)$, defined by $H(l_n(T))=(\lambda_n^X)_*(T)$ for $T\in\KK(X_n)$, where $l_n:\KK(X_n)\longrightarrow\varinjlim\KK(X_n)$  is the canonical map.

Note that $\{\theta_{r,s}: r,s\in   \bigcup\lambda_n^X(X_n)\}$ is dense in  $\KK(X_\infty)$  because $X_\infty=\overline{\bigcup\lambda_n^X(X_n)}$. It follows from that fact that $H$ is onto, since $\theta_{\lambda^X_n(x),\lambda^X_n(y)}=H(l_n(\theta_{x,y}))$, for $x,y\in X_n$.

The map $H$ is also isometric: let $T\in \KK(X_n)$, $T=\sum_{i=1}^k\theta_{x_i,y_i}$, where $x_i,y_i\in X_n$. Then 
\[\|H(l_n(T))\|=\|(\lambda_n^X)_*(T)\|=\|{\textstyle{\sum_1^k}} \theta_{\lambda^X_n(x_i), \lambda^X_n(y_i)}\|=\|(\lambda^A_n)^{(k)}(X^{1/2}Y^{1/2})\|,\]
where (\cite[2.1]{pin}) $X_{ij}=\langle x_i,x_j\rangle$ and  $Y_{ij}= \langle y_i,y_j\rangle$.

Therefore, by applying \cite[2.1]{pin} again,
\[\begin{array}{ll}
\|H(l_n(T))\|&=\|(\lambda^A_n)^{(k)}(X^{1/2}Y^{1/2})\|=\lim_m\|(\phi_{n,m}^A)^{(k)}(X^{1/2}Y^{1/2})\|\\
&=\lim_m\|(\phi_{n,m}^X)_*(T)\|=\|l_n(T)\|.
\end{array}\]
\end{proof}

\section{Cuntz-Pimsner \cstar-algebras and crossed-products by Hilbert \cstar-bimodules}

In this section we show that the pair $(A_\infty,X_\infty)$ obtained in Example \ref{main} is such that $A_\infty\rtimes X_\infty$ is isomorphic to $\ox$. We begin by recalling some well-known facts about adjointable operators on the direct sum of Hilbert \cstar-modules.

Given a sequence $\{X_n\}$ of right Hilbert \cstar-modules over a \cstar-algebra $A$, let $E=\bigoplus_0^\infty X_n$. If $K_0,K_1\subset \NN$, we identify $\LL({\scriptstyle{ \bigoplus_{n\in K_0}X_n, \bigoplus_{n\in K_1}X_n}})$ with a subspace of $\LL(E)$ by extending $\tilde{T}\in \LL({\scriptstyle{ \bigoplus_{n\in K_0}X_n, \bigoplus_{n\in K_1}X_n}})$ to $T\in \LL(E)$ so that $T|_{X_n}=0$ for $n\not\in K_0$. 

Let $J=\overline{\bigcup_m \LL\big( \bigoplus_0^m X_n\big)}\subset \LL(E)$, and let $M$ denote the idealizer of $J$ in $\LL(E)$, that is,  $M=\{T\in\LL(E):TS,ST\in J\ \mbox{for all } S\in J\}$.

For an integer $k$, let 
\[\Delta_k=
\{T\in \LL(E):T(X_n)\subset X_{n+k}\ \mbox{ if } n\geq \max\{0,-k\}, \ T|_{X_n}=0 \mbox{ otherwise}\} .\]

Given  $T\in \Delta_k$, we denote by  $T_n$   the map $T_n\in\LL(X_n,X_{n+k})$ obtained by restricting $T$ to $X_n$. Then $T=\bigoplus_{0}^{\infty}T_n$ and $\|T\|=\sup_n\|T_n\|$.
 Note that $\Delta_k\subset M$ for all $k\in \ZZ$.

\begin{lemma}

\label{qnorm}
If $T\in \Delta_k$, then $\|T+J\|_{M/J}=\limsup_n\|T_n\|$. 

\end{lemma}

\begin{proof}
 We can assume that $k=0$, since $T^*T\in \Delta_0$ and $(T^*T)_n=(T_n)^*T_n$ for all $T\in \Delta_k$. Let $L$ denote $\limsup_n\|T_n\|$.
Given $\epsilon >0$, let $n_0$ be such that $\|T_n\|<L+\epsilon$ for all $n\geq n_0$. 

Then $\|T+J\|_{M/J}\leq \|\bigoplus_{n_0}^\infty T_n\|=\sup_{n\geq n_0}\|T_n\|\leq L+\epsilon,$
which shows that $\|T+J\|_{M/J}\leq L$. On the other hand, if $l<L$ and $S\in \LL(\bigoplus_0^mX_n\big)\subset J$, then
\[\|T-S\|=\|[\bigoplus_0^m T_n-S] \oplus \bigoplus_{m+1}^{\infty}T_n\|\geq \|\bigoplus_{m+1}^{\infty}T_n\|= \sup_n \{\|T_n\|:n>m\|\})>l.\]
Therefore $\|T+J\|_{M/J}\geq l$ for all $l<L$, which ends the proof.

\end{proof}

We next recall the definitions of the Cuntz-Pimsner algebras ${\mathcal{O}}_X$ and $\tilde{{\mathcal{O}}}_X$ given in \cite{pim}.
Given a correspondence $X$ over a \cstar-algebra $A$, let $X_n=X^{\otimes n}$, where $X^{\otimes 0}=A$,  and let $E=\bigoplus_0^{\infty} X_n$ . 

If $x\in X^{\otimes k}$,  we denote by $T_x$  the map $T_x\in\Delta_k\subset \LL(E)$ given by $T_x(y)=x\otimes y$ if $k>0$ and by  $T_a(y)=ay$, if $a\in A$, where $x\otimes a$ is identified with $xa$, for $x\in X^{\otimes k}$, $k\geq 0$,  and $a\in A$.

For $M$ and $J$  defined as above, let $\pi:M\longrightarrow M/J$ be the canonical projection and set $S_x=\pi(T_x)$, for $x\in X_k$, $k\geq 0$. 
The Cuntz-Pimsner \cstar-algebra ${\mathcal{O}}_X$ and the augmented Cuntz-Pimsner \cstar-algebra $\tilde{{\mathcal{O}}}_X$ are the \cstar-subalgebras of $M/J$ generated by $\{S_x:x\in X\}$ and by  $\{S_t:t\in X\cup A\}$, respectively. 
  Notice that $S_{x_1\otimes x_2\otimes \dots\otimes x_k}=S_{x_1} S_{x_2} \dots S_{x_k}$, which implies that $S_x\in {\mathcal{O}}_X$ for all $x\in X^{\otimes k}$, $k\geq 1$.

\begin{rk}
\label{snorm}
{\rm{Let $x\in X^{\otimes m}$, for $m\geq 0$. Since $\|(T_x)_n\|=\|(T_x)_{n-1}\otimes \id_X\|\leq \|(T_x)_{n-1}\|$ for all $n\geq 1$, we have  by Lemma \ref{qnorm}
\[\|S_x\|=\limsup_n\|(T_x)_{n}\|=\inf_{n\geq 1}\|(T_x)_{k}\otimes \id_{X^{\otimes n}}\|,\]
  for all  $k\geq 1$.}}
\end{rk}

\begin{lemma}
Let $(X_i,\phi_i)$ be $A-B_i$ correspondences for $i=1,2$, and let $Y$ be  a right Hilbert \cstar-module over  $A$. If $\ker\phi_1\subset \ker\phi_2$, then $\|T\otimes \id_{X_1}\|\geq \|T\otimes \id_{X_2}\|$ for all $T\in \LL(Y)$.

\end{lemma}
\begin{proof}
Let $T\in {\mathcal{L}}(Y)$.  Then $T\otimes \id_{X_i}=0$ if and only if $0=\|Ty\otimes x\|^2=\|\langle x,\langle Ty,Ty\rangle \cdot x\rangle\|$ for all $x\in X_i, \ y\in Y$. That is, $T\otimes \id_{X_i}=0$ if and only if $\langle Ty,Ty\rangle\in \ker \phi_i$ for all $y\in Y$. We can thus define a map $T\otimes \id_{X_1}\mapsto T\otimes \id_{X_2}$, which is a (norm-decreasing) *-homomorphism between the \cstar-algebras $\{T\otimes \id_{X_1}: T\in {\mathcal{L}}(Y)\}$ and $\{T\otimes \id_{X_2}: T\in {\mathcal{L}}(Y)\}$. 
\end{proof}

\begin{clly}
\label{id}
Let $(X,\phi_X)$ and $Y$ be, respectively,  a correspondence and a right Hilbert \cstar-module  over a \cstar-algebra $A$. Let $X(A)$ be as in Example \ref{main}.  Then for any $T\in {\mathcal{L}}(Y)$ we have 
\[ \|T\otimes \id_X\|=\|T\otimes \id_{\rxa}\|.\]
\end{clly}
\begin{proof}
It suffices to notice that an element $a$ of $A$ acts on $\rxa$ by left multiplication by $\phi_X(a)$. Since Im $\phi_X\subset \rxa$, we conclude that $a\cdot\rxa=0$ if and only if $\phi_X(a)=0$. Then the previous lemma applies in both directions and the equality holds.

\end{proof}

\begin{lemma} 
\label{ad} Given  a correspondence  $(X,\phi_X)$ over a \cstar-algebra $A$,  let  $X(X)$,  $X(A)$, and  $\psi_X:X\longrightarrow X(X)$  be as in Example \ref{main}.

For $n\geq 1$, let  $\beta_n:\rx^{\otimes n}\longrightarrow X^{\otimes n}\otimes_{\phi_X} \rxa$ be the isomorphism of $\rxa$- correspondences  given by $\beta_n=\id_X\otimes L_{X(A),X}^{\otimes n-1}\otimes \id_{\rxa}$, where  $L_{X(A),X}$ is as in Notation \ref{not}.  
 Then:
\begin{enumerate}

\item $\beta_n\circ \psi_X^{\otimes n}=\id_{X^{\otimes n-1}}\otimes \psi_X$, for all $n\geq 1$.

\item   $\beta_{n+m}\big(\theta_{\psi_X^{\otimes n}(z), \psi_X^{\otimes n}(w)}\otimes \id_{\rx^{\otimes m}}\big)\beta^*_{n+m}=\theta_{z,w}\otimes \id_{(X^{\otimes m}\otimes \rxa)}$,
for $z,w\in X^{\otimes n}$, $n\geq 1$, $m\geq 0$.

\item  $\beta_{m+1}\big(\phi_X(a)\otimes \id_{\rxa\otimes \rx^{\otimes m}}\big)\beta^*_{m+1}=\phi_X(a)\otimes \id_{(X^{\otimes m}\otimes \rxa)}$, for all $m\geq 0$.

\end{enumerate}

\end{lemma}

\begin{proof}

(1) Let $x_i\in X$, $a_i\in A$, for $i=1,2,\dots ,n$. Then

\[[\beta_n\circ \psi_X^{\otimes n}](x_1a_1\otimes x_2a_2\otimes\dots\otimes x_na_n)\]
\[=\beta_n(x_1\otimes \phi_X(a_1)\otimes  x_2\otimes \phi_X(a_2)\otimes \dots \otimes x_n\otimes \phi_X(a_n))\]
\[=x_1\otimes \phi_X(a_1)  x_2\otimes \phi_X(a_2)x_3\otimes \dots \otimes \phi_X(a_{n-1})x_n\otimes \phi_X(a_n)\]
\[= x_1a_1\otimes x_2a_2\otimes\dots\otimes x_{n-1}a_{n-1}\otimes x_n\otimes \phi_X(a_n)\]
\[= (\id_{X^{\otimes n-1}}\otimes \psi_X)(x_1a_1\otimes x_2a_2\otimes\dots\otimes x_na_n).\]

(2) We first prove the statement for  $m=0$. We assume, without loss of generality, that $z=z_0\otimes xa$, $w=w_0b$, for $x\in X$, $z_0\in X^{\otimes n-1}$,  $w_0\in X^{\otimes n}$, and $a,b\in A$.
Let $u\in X^{\otimes n}$, $r\in\rxa$. Then, by (1):

\[\big(\beta_{n}\theta_{\psi_X^{\otimes n}(z), \psi_X^{\otimes n}(w)}\beta^*_{n}\big)(u\otimes r)\]
\[= \big(\theta_{(\id_{X^{\otimes n-1}}\otimes\psi_X)(z_0\otimes xa),(\id_{X^{\otimes n-1}}\otimes \psi_X) (w_0b)}\big)(u\otimes r)\]
\[=z_0\otimes x\otimes \phi_X(a)\langle w_0\otimes \phi_X(b),u\otimes r\rangle\]
\[=z_0\otimes xa\otimes \phi_X(b^* \langle w_0,u\rangle) r\]
\[=z_0\otimes xa\langle w_0b,u\rangle\otimes r=z\langle w,u\rangle \otimes r\]
\[=(\theta_{z,w}\otimes \id_{\rxa})(u\otimes r).\]

Let us denote by $L$ the map $ L_{X(A),X}$ defined in Notation \ref{not}. For $m\geq 1$ we have
\[\beta_{n+m}=\big(\id_{X^{\otimes n}}\otimes L^{\otimes m}\otimes \id _{\rxa}\big)\big(\beta_n\otimes \id_{\rx^{\otimes m}}\big).\]

Therefore 
\[\beta_{n+m}(\theta_{\psi^{\otimes n}_X(z),\psi^{\otimes n}_X(w)}\otimes \id_{\rx^{\otimes m}}\big)\beta^*_{n+m}=\]
\[\big(\id_{X^{\otimes n}}\otimes L^{\otimes m}\otimes \id _{\rxa}\big)(\theta_{z,w}\otimes \id_{X(A)}\otimes \id_{X(X)^{\otimes m}})\big(\id_{X^{\otimes n}}\otimes (L^{\otimes m})^*\otimes \id_{\rxa}\big)\]
\[=\theta_{z,w}\otimes \id_{X^{\otimes m}\otimes X(A)}.\]

(3) 
\[\beta_{m+1}\big(\phi_X(a)\otimes \id_{\rxa\otimes \rx^{\otimes m}}\big)\beta^*_{m+1}=\]
\[(\id_X\otimes L^{\otimes m}\otimes \id_{\rxa})(\phi_X(a)\otimes \id_{(\rxa\otimes X)^{\otimes m}}\otimes \id_{\rxa})(\id_X\otimes(L^{\otimes m})^*\otimes \id_{\rxa})\]
\[=\phi_X(a)\otimes \id_{X^{\otimes m}\otimes \rxa}.\]

\end{proof}

\begin{rk}
\label{dual}

{\rm{As discussed in \cite[Remark 1.2, (2)]{pim},  the automorphism of $X$ given  by $x\mapsto\lambda x$ for $\lambda\in S^1$ yields an automorphism $\gamma_\lambda$ of $\tilde{{\mathcal{O}}}_X$, determined by $\gamma_\lambda(S_x)=\lambda^kS_x$ for $x\in X^{\otimes k}$, $k\geq 0$. In fact, this automorphism of $X$  extends to an automorphism $\dot{\gamma}_\lambda$ of $E$ defined by $\big(\dot{\gamma}_\lambda(\eta)\big)(k)=\lambda^k\eta(k)$, for $\eta\in E$. Conjugation by $\dot{\gamma}_\lambda$ is an automorphism of $\LL(E)$ that maps $T_x$ into $\lambda^kT_x$ for $x\in X^{\otimes k}$, $k\geq 0$, and it  leaves $J$ invariant.

Thus one gets an action $\gamma$ of $S^1$ on $\tilde{{\mathcal{O}}}_X$ that is easily checked to be strongly continuous. The fixed-point subalgebra of this action is $E_0(\ox)=\overline{\rm{span}}\{S_xS^*_y: x,y\in X^{\otimes n},\ n\geq 0\}$ and its first spectral subspace $E_1(\ox)=\overline{\rm{span}}\{S_xS^*_y: x\in X^{\otimes n+1},\ y\in X^{\otimes n},\ n\geq 0\}=\overline{\rm{span}}\{S_xe:x\in X,\  e\in E_0(\ox)\}$.
  This last statement is shown by means of the usual argument, since span$\{S_xS^*_y: x\in X^{\otimes n}, y\in X^{\otimes m},\ n,m\geq 0\}$ is dense in $\ox$, and the maps $P_i:\ox\rightarrow E_i(\ox)$ given by 
\[P_i(u)=\int_{S^1} z^{-i}\gamma_z(u)dz \mbox{ for }u\in\ox\]
are surjective contractions (see \cite{circ} for details), and $\gamma_{\lambda}(S_xS^*_y)=\lambda^{n-m}S_xS^*_y$, for $x\in X^{\otimes n} $,  $y\in X^{\otimes m} $, and $n,m \geq 0$, $i=0,1$.

Now, since  $\tilde{{\mathcal{O}}}_X$ is generated as a \cstar-algebra by $E_0(\ox)$ and $E_1(\ox)$,  Theorem 3.1 in \cite{aee} applies, and $\tilde{{\mathcal{O}}}_X$ is isomorphic to the crossed-product $E_0(\ox)\rtimes E_1(\ox)$.}}

\end{rk}
\begin{prop}
\label{isoo}
Let $(X, \phi_X)$ be a correspondence over a \cstar-algebra $A$, and let $X(A)$, $X(X)$  and $\psi_X$ be as in Example \ref{main}. Then there is an isomorphism of Hilbert \cstar-bimodules $(\eta_1,\eta_0):\big(E_1(\tilde{{\mathcal{O}}}_X), E_0(\tilde{{\mathcal{O}}}_X)\big)\longrightarrow \big(E_1(\tilde{{\mathcal{O}}}_{X(X)}), E_0(\tilde{{\mathcal{O}}}_{X(X)})\big)$  carrying $S_x$ and  $S_a$ to  $S_{\psi_X(x)}$ and $S_{\phi_X(a)}$  respectively, for $x\in X$ and $a\in A$.

Besides,   if  $(i_X,i_A):(X,A)\longrightarrow \tilde{{\mathcal{O}}}_X$ is  given by $i_X(x)=S_x$ and  $i_A(a)=S_a$ and similarly for $(X(X),X(A))$, then
\[\xymatrix{
{\scriptstyle{(X,A)}} \ar[d]_{(i_X,i_A)}\ar[r]^{(\psi_X,\phi_X)\ \ \ \ } &{\scriptstyle{  (X(X),X(A))}}\ar[d]^{(i_{X(X)},i_{X(A)})}\\
{\scriptstyle{\big(E_1(\tilde{{\mathcal{O}}}_X), E_0(\tilde{{\mathcal{O}}}_X)\big)}}\ar[r]^{(\eta_1,\eta_0)} & {\scriptstyle{\big(E_1(\tilde{{\mathcal{O}}}_{X(X)}), E_0(\tilde{{\mathcal{O}}}_{X(X)})\big)}} }\]
is a commuting diagram of homomorphisms of correspondences.
\end{prop}

\begin{proof}

We would like to define
$\eta_0:E_0(\ox)\longrightarrow E_0({\mathcal{{\tilde{O}}}}_{X(X)})$ by 
\[\eta_0(S_a+\sum_{i=1}^kS_{x_i}S^*_{y_i}):= S_{\phi_X(a)}+\sum_{i=1}^kS_{\psi_X^{\otimes n_i}(x_i)}S^*_{\psi_X^{\otimes n_i}(y_i)},\]
 where $a\in A$, $x_i,y_i\in X^{\otimes n_i}$, and  $n_i>0$ for all $i=1,...,k$.

We first show that the definition above makes sense. Let $a,x_i,y_i$ be as above, and let $m=\max\{n_i:i=1,...,k\}$. Then (see the beginning of Section 1 in  \cite{pim} for the first equality)

\[\begin{array}{ll}
T_a+\sum_{i=1}^kT_{x_i}T^*_{y_i}&=\bigoplus_{n=m}^\infty (T_a)_n+\big(\sum_i \theta_{x_i,y_i}\otimes \id_{X^{\otimes n-n_i}}\big)\ \ \mbox{{\em (mod $J$)}}\\
&=\bigoplus_{n=m}^\infty \tau(\{a,x_i,y_i\}) \otimes \id_{X^{\otimes n-m}},
\end{array}\]
where $\tau(\{a,x_i,y_i\})=\phi_X(a)\otimes\id_{X^{\otimes m-1}}+\sum_i \theta_{x_i,y_i}\otimes \id_{X^{\otimes m-n_i}}$ .

Now, by parts (2) and (3) of Lemma \ref{ad}:
\[\begin{array}{c}
\tau(\{a,x_i,y_i\})\otimes \id_{X^{\otimes n}\otimes X(A)}\\
=\phi_X(a)\otimes \id_{X^{\otimes m-1+n}}\otimes \id_{X(A)}+\sum_i\theta_{x_i,y_i}\otimes \id_{X^{\otimes m-n_i+n}}\otimes\id_{X(A)}\\
=\beta_{n+m}\big(\phi_X(a)\otimes\id_{X(A)}\otimes \id_{X(X)^{\otimes m-1+n}} +\\
+\sum_i\theta_{\psi_X^{\otimes n_i}(x_i),\psi_X^{\otimes n_i}(y_i)}\otimes \id_{(X(X))^{\otimes m-n_i+n}}\big)\beta^*_{n+m}\\
=\beta_{n+m}\big(\tau(\{\phi_X(a),\psi_X^{\otimes n_i}(x_i), \psi_X^{\otimes n_i}(y_i)\}\big)\otimes \id_{X(X)^{\otimes n}})\beta^*_{n+m}.
\end{array}\]
Therefore, by Lemma  \ref{qnorm} and Corollary \ref{id}
\[\begin{array}{ll}
\|S_a+\sum_iS_{x_i}S^*_{y_i}\|&=\lim_n \|\tau(\{a,x_i,y_i\})\otimes \id_{X^{\otimes n}}\|\\
&=\lim_n \|\tau(\{a,x_i,y_i\})\otimes \id_{X^{\otimes n}}\otimes \id_{X(A)}\|\\
&=\lim_n \|\tau(\{\phi_X(a),\psi_X^{\otimes n_i}(x_i), \psi_X^{\otimes n_i}(y_i)\})\otimes \id_{X(X)^{\otimes n}}\|\\
&=\|\eta_0\big(S_a+\sum_iS_{x_i}S^*_{y_i}\big)\|.
\end{array}\]
This shows that $\eta_0$ can be extended to an isometry $\eta_0:E_0(\ox)\longrightarrow E_0({\mathcal{{\tilde{O}}}}_{X(X)})$  which is  easily checked to be an isometric *-homomorphism, in view of the properties listed in \cite[Proposition 1.3]{pim}.

We next show that $\eta_0$ is onto. First note that $S_{\psi_X(x)}S^*_{\psi_X(y)}=S_{\theta_{x,y}}$ for all $x,y\in X$. Let $\pi:M\rightarrow M/J\supset {\mathcal{{\tilde{O}}}}_{X(X)}$ be as in the beginning of this section. 
By \cite[1.3]{pim}:
\[S_{\psi_X(x)}S^*_{\psi_X(y)}=\pi\big(\oplus_{n=0}^\infty\theta_{\psi_X(x),\psi_X(y)}\otimes \id_{X(X)^{\otimes n}}\big)=S_{\theta_{x,y}}\]
because
\[\theta_{\psi_X(x),\psi_X(y)}=\theta_{x,y}\otimes\id_{X(A)}=\phi_{X(X)}\big(\theta_{x,y}\big) .\]
Also,  by \cite[1.3]{pim}
\[S_{xa\otimes \theta_{y,z}}=S_{[x\otimes \phi(a)]\theta_{y,z}}= S_{\psi_X(xa)}S_{\theta_{y,z}}=S_{\psi_X(xa)}S_{\psi_X(y)}S^*_{\psi_X(z)}.\]

Since $S_{\phi_X(a)}=\eta_0(S_a)$, $S_{\theta_{x,y}}=\eta_0(S_xS^*_y)$,  and $S_{u_1\otimes u_2\otimes \dots \otimes u_n}=S_{u_1}S_{u_2}\dots S_{u_n}$, it only remains to show that $S_{u_1}S_{u_2}\dots S_{u_n}S^*_{v_n}\dots S^*_{v_2} S^*_{v_1}\in {\rm Im}\ \eta_0$ for all   $u_i, v_i\in X(X)$ and $n\geq 1$. We proceed by induction on $n$.
The case $n=1$ follows from the fact that, by the identities above:
\[\begin{array}{c}
S_{x\otimes (\phi_X(a)+\theta_{y,z})}S^*_{x'\otimes (\phi_X(a')+\theta_{y',z'})}\\
= \big( S_{\psi_X(xa)}+S_{\psi_X(x)}S_{\psi_X(y)}S^*_{\psi_X(z)}\big)\big( S_{\psi_X(x'a')}+S_{\psi_X(x')}S_{\psi_X(y')}S^*_{\psi_X(z')}\big)^*\\
=\eta_0\big(S_{xa}S^*_{x'a'}+S_{xa}S_{z'}S_{y'}^*S_{x'}^*+S_{x}S_{y}S^*_{z}S^*_{x'a'}+ S_{x}S_{y\langle z,z'\rangle}S_{ y'}^*S_{x'}^*\big).
\end{array}\]
for all $a,a'\in A$ and $x,x',y,y',z,z'\in X$.
The induction step follows from the fact  that for all $a,a'\in A$ and $x,x',y,y',z,z'\in X$
\[S_{x\otimes (\phi_X(a)+\theta_{y,z})}\big({\rm Im }\  \eta_0\big) S^*_{x'\otimes (\phi_X(a')+\theta_{y',z'})}\subset {\rm Im }\  \eta_0,\]
which is checked by applying the action $\gamma$ of Remark \ref{dual} (or by direct computation).

We now define  
\[\eta_1:E_1(\ox)\longrightarrow E_1({\mathcal{{\tilde{O}}}}_{X(X)}) \mbox{ by }\eta_1(\sum_iS_{x_i}e_i)=\sum_i S_{\psi_X(x_i)}\eta_0(e_i),\]
 for $x_i\in X$ and $e_i\in E_0(\ox)$. To check that the map $\eta_1$ thus defined makes sense and extends to an isometric map on $E_1(\ox)$ notice that
\[\begin{array}{ll}
\|\sum_iS_{\psi_X(x_i)}\eta_0(e_i)\|^2&=\|\sum_{i,j}\eta_0(e_i)^* S_{\psi_X(x_i)}^*S_{\psi_X(x_j)}\eta_0(e_j)\|\\
&=\|\sum_{i,j}\eta_0\big(e_i^*S_{x_i}^*S_{x_j}e_j\big)\|\\
&=\|\sum_iS_{x_i}e_i\|^2.
\end{array}\] 
Straightforward computations show that $(\eta_0, \eta_1)$ is a Hilbert \cstar-bimodule homomorphism.  Besides, the map $\eta_1$ is  onto because so is $\eta_0$. It is clear from the definitions that the diagram commutes. Finally, it follows from \cite[Proposition 1.3]{pim} that $(i_X,i_A)$ is  a homomorphism of correspondences and it was shown in Example \ref{main} that so is $(\psi_X,\phi_X)$. 
\end{proof}

\begin{clly} Let $(X,\phi_X)$ be a correspondence over a \cstar-algebra $A$,  and let $X(A)$ and $X(X)$ be as in Example \ref{main}. Then $\ox$ and ${\mathcal{\tilde{O}}}_{X(X)}$ are isomorphic.
\end{clly}
\begin{proof}
The isomorphism of Hilbert \cstar-bimodules $(\eta_1,\eta_0)$ obtained in Proposition \ref{isoo} induces an isomorphism from $E_0(\ox)\rtimes E_1(\ox)$ to $E_0({\mathcal{\tilde{O}}}_{X(X)})\rtimes E_1({\mathcal{\tilde{O}}}_{X(X)})$. The statement now follows from Remark \ref{dual}.
\end{proof}

\begin{thm}
\label{eqcp}
Let $X$ be a correspondence over a \cstar-algebra $A$, and let $(X_\infty,A^X_\infty)$ be as in Example \ref{main}. Then $\tilde{{\mathcal{O}}}_X\cong A^X_{\infty}\rtimes X_{\infty}$.

\end{thm}
\begin{proof}
As in Example \ref{main}, let   $(X_\infty,A^X_\infty,\{\lambda_n^X\}, \{\lambda_n^A\})$  be the direct limit of the directed sequence  $\{(X_n,A^X_n,\phi_n^X,\phi_n^{A,X})\}$, and let 
\[( i_{X_n}, i_{A^X_n}):(X_n,A^X_n)\rightarrow \big(E_1(\tilde{{\mathcal{O}}}_{X_n}), E_0(\tilde{{\mathcal{O}}}_{X_n})\big)\] and \[(\eta^n_1,\eta^n_0):\big(E_1(\tilde{{\mathcal{O}}}_{X_n}), E_0(\tilde{{\mathcal{O}}}_{X_n})\big)\longrightarrow \big(E_1(\tilde{{\mathcal{O}}}_{X_{n+1}}), E_0(\tilde{{\mathcal{O}}}_{X_{n+1}})\big)\]  be as in Proposition \ref{isoo}. 

Let  $\Upsilon^n_i:E_i(\tilde{{\mathcal{O}}}_{X_n}) \rightarrow E_i(\tilde{{\mathcal{O}}}_X)$ be the isomorphism of Hilbert \cstar-bimodules defined by  $\Upsilon^n_i=(\eta^0_i\eta^1_i\dots\eta^{n-1}_i)^{-1}$ for all $n\geq 0$, $i=0,1$.

By  Propositions \ref{isoo} and \ref{dlim} and Remark \ref{bimodule} there are homomorphisms of Hilbert \cstar-bimodules $(i_\infty^X, i_\infty^A)$  making the diagram
\[
\xymatrix{
(X_n,A^X_n)\ar[dr]^{(\lambda_n^X, \lambda_n^A)} \ar[rr]^{(\phi_n^X,\phi_n^{A,X})\ \ \ }\ar[dddr]_{(\Upsilon^n_1i_{X_n},\Upsilon^n_0i_{A^X_n}) } & & (X_{n+1}, A^X_{n+1}) \ar[dl]_{(\lambda_{n+1}^X, \lambda_{n+1}^A)}\ar[dddl]^{(\Upsilon^{n+1}_1i_{X_{n+1}},\Upsilon^{n+1}_0i_{A^X_{n+1}})}\\
& (X_{\infty},A^X_{\infty})\ar@{-->}[dd]|{(i_\infty^X,i_\infty^A)} & \\
& &\\
& {\scriptstyle{\big(E_1(\tilde{{\mathcal{O}}}_X), E_0(\tilde{{\mathcal{O}}}_X)\big)}}  & }
\]
commute.

Since by Proposition \ref{dlim} the  pair  $(i_\infty^X, i_\infty^A):(X_\infty,A_\infty^X)\longrightarrow \tilde{{\mathcal{O}}}_X$ is  covariant in the sense of \cite[2.1]{aee}, it induces, by the universal property of the crossed product, a *-homomorphism $i: A^X_\infty\rtimes X_\infty\longrightarrow \tilde{{\mathcal{O}}}_X$, which  is onto because its image contains $\{S_x:x\in X\cup A\}$.
It only remains to check that $i^A_\infty$ is injective, since this would imply by \cite[2.9]{circ} that so is $i$, $i$ being covariant for the dual action (\cite[3]{aee}) on the crossed product and the action $\gamma$ discussed in Remark \ref{dual}  on $\tilde{{\mathcal{O}}}_X$.

First notice that 
\[\|\phi^X_n(a_n)\otimes \id_{X_n^{\otimes k}}\|=\|a_n\otimes \id_{X_{n-1}^{\otimes k+1}}\|,\]
for $n\geq 1,\ k\geq 0,\ a_n\in A^X_n\subset \LL(X_{n-1})$. In fact, the unitary $\id_{X_{n-1}}\otimes L^{\otimes k}_{A_n^X,X_{n-1}}\otimes \id_{A^X_n}$, for $L_{A_n^X,X_{n-1}}$ as in Notation \ref{not}, intertwines $\phi^X_n(a_n)\otimes \id_{X_n^{\otimes k}}$ and $a_n\otimes \id_{X_{n-1}^{\otimes k}}\otimes \id_{A_n}$. Now, by Corollary \ref{id}:
\[\|\phi^X_n(a_n)\otimes \id_{X_n^{\otimes k}}\|=\|a_n\otimes \id_{X_{n-1}^{\otimes k}}\otimes \id_{A_n^X}\|=\|a_n\otimes \id_{X_{n-1}^{\otimes k+1}}\|.\]
It now follows by induction on $m-n$ that
\[\|\phi^X_{n,m}(a_n)\|=\|a_n\otimes \id_{X_{n-1}^{\otimes m-n}}\|,\]
for $m\geq n\geq 1$ and  $ a_n\in A_n^X$.

We next show that $i_\infty^A$ is injective by showing that its restriction to $\lambda_n^A(A^X_n)$ is isometric for all $n\geq 1$.
Take $a_n\in A_n^X$ for $n\geq 1$. Then:
\[\begin{array}{ll}
\|\lambda^A_n(a_n)\|&=\lim_m\|\phi^X_{n,m}(a_n)\|\\
&=\lim_m\|a_n\otimes \id_{X_{n-1}^{\otimes m-n}}\|\\
&=\lim_m\|a_n\otimes \id_{X_{n-1}^{\otimes m}}\|\\
&=\|S_{a_n}\|\\
&=\|i_{A^X_n}(a_n)\|\\
&=\|i^A_{\infty}(\lambda^A_n(a_n))\|
\end{array}\]
\end{proof}
\section{Morita equivalence for Cuntz-Pimsner \cstar-algebras}

We establish in this section a sufficient condition for the Morita equivalence of two augmented Cuntz-Pimsner \cstar-algebras. In order to do so,  we view these algebras as crossed products by Hilbert \cstar-bimodules as in Theorem \ref{eqcp}, and then we use the condition for the Morita equivalence of crossed products given in \cite[4,2]{aee}.
Along this section we will be making extensive use of the construction described in  Example \ref{main}.

\begin{lemma}
\label{pimsner}
Let $(Y,\phi_Y)$ be a correspondence over a  \cstar-algebra $B$ and let $M$ be a right Hilbert \cstar-module over $B$ . For  $Y(B)$ and $Y(M)$  as in Example \ref{main}, there is a $*$-homomorphism $O_1:\KK(M\otimes Y)\longrightarrow \KK(Y(M))$ such that
\[[O_1(\theta_{m_1\otimes y_1, m_2\otimes y_2})](m\otimes r)=m_1\otimes \theta_{y_1, y_2} \phi_Y(\langle m_2,m\rangle )r,\]
for all $m, m_1,m_2\in M$, $y_1,y_2\in Y$, and $r\in Y(B)$, where  $Y(M)$ is  viewed as a $Y(B)$-right Hilbert \cstar module.

Besides,  $O_1\big(\theta_{m_1\otimes y_1\langle y_2,z_2\rangle, m_2\otimes z_1}\big)=\theta_{m_1\otimes\theta_{y_1,y_2}, m_2\otimes\theta_{z_1,z_2}}.$

\end{lemma}
\begin{proof}

It was shown in \cite[2.2]{pim} that $\KK(M\otimes Y)$ and $\KK(M\otimes \KK(Y))$ are isomorphic. Now  $\KK(M\otimes \KK(Y))$  can be viewed as contained in $\KK(Y(M))$, since $M\otimes \KK(Y)$ is a closed  $Y(B)$-right Hilbert \cstar-submodule  of $ Y(M)$: In fact, if $x_i,y_i\in M\otimes \KK(Y)$ for $i=1,\dots,n$, we have by \cite[2.1]{pin}:
\[\|\sum_{i=1}^n\theta^{Y(M)}_{x_i,y_i}\|=\|A^{1/2}C^{1/2}\|_{M_n(Y(B))}= \|\sum\theta^{M\otimes \KK(Y)}_{x_i,y_i}\|,\]
where $A_{ij}= \langle x_i,x_j\rangle^{Y(M)}=\langle x_i,x_j\rangle^{M\otimes \KK(Y)}$, and anagously for $C$. 

In this way we can obtain an isometric *-homomorphism  $I:\KK(M\otimes \KK(Y))\hookrightarrow \KK(Y(M))$, defined by $I(\theta^{Y(M)}_{x,y})=\theta^{M\otimes \KK(Y)}_{x,y}$. The map $O_1$ is now defined to be the composition of  the isomorphism $P$ in \cite[2.2]{pim} with $I$.  By  keeping track of the proof in \cite[2.2]{pim}, we get the formulas in the statement.  In fact, let us identify  $w_1\otimes \tilde{w}_2$ with $\theta_{w_1,w_2}$, for $w_i\in Y$, $i=1,2$. 

Then, according to \cite[2.2]{pim},  $[P\big(\theta_{m_1\otimes y_1, m_2\otimes y_2}\big)](m\otimes \theta_{w_1,w_2})]$ gets identified with
\[(\theta_{m_1\otimes y_1, m_2\otimes y_2}(m\otimes w_1))\otimes \tilde{w}_2=m_1\otimes y_1\langle y_2, \phi_Y(\langle m_2,m\rangle) w_1\rangle \otimes \tilde{w}_2,\]
which gets identified with $m_1\otimes \theta_{y_1,y_2}\phi_Y(\langle m_2,m\rangle)\theta_{w_1,w_2}$.

Straightforward computations now show that  
\[[P\big(\theta_{m_1\otimes y_1\langle y_2,z_2\rangle, m_2\otimes z_1}\big)](\xi)=\theta_{m_1\otimes\theta_{y_1,y_2}, m_2\otimes\theta_{z_1,z_2}}(\xi)\]
when $\xi\in M\otimes \KK(Y)$. Then, by applying the map $I$ we get: 
\[O_1\big(\theta_{m_1\otimes y_1\langle y_2,z_2\rangle, m_2\otimes z_1}\big)=\theta_{m_1\otimes\theta_{y_1,y_2}, m_2\otimes\theta_{z_1,z_2}},\] 
which yields the formulas in the statement.

\end{proof}

\begin{prop}
\label{bigO}
Let $(Y,\phi_Y)$ be a correspondence over a  \cstar-algebra $B$, and let $M$ be a right Hilbert \cstar-module over $B$.
Let $L_1$ and $L_2$ be the \cstar-subalgebras of  $\LL(M\otimes Y)$ defined  by $L_1=\KK(M\otimes Y)$ and $L_2= \{T\otimes \id_Y: T\in \KK(M)\}$, and  let $L=L_1+L_2$ be  the \cstar-subalgebra  of  $\LL(M\otimes Y)$  generated by $L_1\cup L_2$.
Then there is  an isomorphism $O:L\longrightarrow \KK(Y(M))$.

\end{prop}

\begin{proof}
We first set $O_i:L_i\rightarrow \KK(Y(M))$, for $i=1,2$ as follows: $O_1$ is the *-homomorphism defined in Lemma \ref{pimsner}  and, in view of Corollary  \ref{id}, we set  $O_2(T\otimes \id_Y)=T\otimes \id_{Y(B)}$,  for $T\in \KK(M)$.

Our aim is to define $O(T_1+T_2)=O_1(T_1)+O_2(T_2)$  for $T_i\in L_i$, $i=1,2$.
To make sense of this, first  note that 
\[O_i(T)\otimes \id_Y=(\id_M\otimes L_{Y(B),Y})^{-1}T (\id_M\otimes L_{Y(B),Y}),\]
for $T\in L_i$, $i=1,2$, and $L_{Y(B),Y}$ as in Notation \ref{not}

The equality is easily checked  for $i=2$ whereas, if $T=\theta_{m_1\otimes y_1,m_2\otimes y_2}$, for   $m_1,m_2\in M$, $y_1,y_2\in Y$,  and  $m\otimes r\otimes y\in M\otimes Y(B)\otimes Y$, then:
\[\begin{array}{l}
[(\id_M\otimes L_{Y(B),Y})^{-1}T (\id_M\otimes L_{Y(B),Y})](m\otimes r\otimes y)\\
=m_1\otimes L^{-1}_{Y(B),Y}\big(y_1\langle m_2\otimes y_2, m\otimes ry\rangle\big)\\
=m_1\otimes L^{-1}_{Y(B),Y}\big(y_1\langle y_2,\phi_Y(\langle m_2,m\rangle)ry\rangle\big)\\
=m_1\otimes \theta_{y_1,y_2}\phi_Y(\langle m_2,m\rangle)r \otimes y\\
=(O_1(T)\otimes \id_Y)(m\otimes r\otimes y).
\end{array}\]
On the other hand, it is straightforwardly verified that 
\[O_i(T)\otimes\id_{Y(B)}=(\id_M\otimes L_{Y(B),Y(B)})^{-1}O_i(T)(\id_M\otimes L_{Y(B),Y(B)}),\]
for $T\in L_i$, $i=1,2$.

By virtue of Corollary \ref{id} and the identities above we  have, for $T_i\in L_i$, $i=1,2$:
\[\begin{array}{c}
\|O_1(T_1)+O_2(T_2)\|=\\
=\| (\id_M\otimes L_{Y(B),Y(B)})[\big(O_1(T_1)+O_2(T_2)\big)\otimes\id_{Y(B)}](\id_M\otimes L_{Y(B),Y(B)})^{-1}\|\\
=\|\big(O_1(T_1)+O_2(T_2)\big)\otimes \id_{Y(B)}\|\\
=\|\big(O_1(T_1)+O_2(T_2)\big)\otimes \id_Y\|\\
=\|(\id_M\otimes L_{Y(B),Y})^{-1}(T_1+T_2) (\id_M\otimes L_{Y(B),Y})\|\\
=\|T_1+T_2\|,
\end{array}\]
which shows that $O$ can be defined as above, and it is an isometric linear map that preserves the involution.

Now, if $T_i\in L_i$, $T_1=\theta_{m_1\otimes y_1, m_2\otimes y_2}$ and $T_2=S\otimes \id_Y$, then
\[\begin{array}{ll}
O_1(T_2T_1)&=O_1(\theta_{Sm_1\otimes y_1, m_2\otimes y_2})\\
&=(S\otimes \id_{Y(B)})O_1(\theta_{m_1\otimes y_1, m_2\otimes y_2})\\
&=O_2(T_2)O_1(T_1).
\end{array}\]

It follows from this and from the fact that $O_1$ and $O_2$ are *-homomorphisms that  $O$ is multiplicative. It only remains to show that $O$ is onto. This fact follows from the following identities that can be verified directly from the definitions:

\begin{itemize} 
\item $\theta_{m_1\otimes \phi_Y(b_1), m_2\otimes \phi_Y(b_2)}=\theta_{m_1b_1,m_2b_2}\otimes \id_{Y(B)}=O\big(\theta_{m_1b_1,m_2b_2}\otimes \id_Y\big)$
\item   $\theta_{m_1\otimes \theta_{y_1,y_2},m_2\otimes \phi_Y(b)}=O\big(\theta_{m_1\otimes y_1,m_2b\otimes y_2}\big)$ 
\item $\theta_{m_1\otimes\theta_{y_1,y_2}, m_2\otimes\theta_{z_1,z_2}}=O\big(\theta_{m_1\otimes y_1\langle y_2,z_2\rangle, m_2\otimes z_1}\big)$,
\end{itemize}
where  $m_1,m_2\in M$, $y_1,y_2, z_1,z_2\in Y$, and $b,b_1,b_2\in B$.
\end{proof}

\begin{rk}
\label{oops}
{\rm{ Notice that we have shown at the beginning of the proof of Proposition \ref{bigO} the identity 
\[(\id_M\otimes L_{Y(B),Y})^{-1}T (\id_M\otimes L_{Y(B),Y})=O(T)\otimes \id_Y\]
for any $T$ belonging to the  \cstar-subalgebra of $\LL(M\otimes Y)$ generated by $\KK(M\otimes Y)\cup \{T\otimes \id_Y: T\in \KK(M)\}$.}}

\end{rk}
\begin{prop}
\label{i1}
Let $X$ and $Y$ be correspondences  over \cstar-algebras $A$ and $B$, respectively, and let $M$ be an $A-B$ Hilbert \cstar-bimodule that is full on the left and  such that there is an isomorphism  $J:X\otimes M\longrightarrow M\otimes Y$ of $A-B$ correspondences.

Let $I:X(A)\longrightarrow \KK(Y(M))$ be given by $I(T)=O(J(T\otimes \id_M)J^{-1})$, for $O$ as in Proposition \ref{bigO}. Then $I$ is an isomorphism, and $I(\phi_X(a))=\phi_M(a)\otimes \id_{Y(B)}$ for all $a\in A$.

\end{prop}

\begin{proof}

By Proposition \ref{bigO}, it suffices to show that $T\mapsto J(T\otimes \id_M)J^{-1}$ is an isomorphism from $X(A)$ to the \cstar-subalgebra $L$ of $\LL(M\otimes Y)$ generated by $\KK(M)\otimes \id_Y$ and $\KK(M\otimes Y)$. 

The image of $X(A)$ by the map $T\mapsto T\otimes\id_M$ is  the \cstar-subalgebra $C$ of $\LL(X\otimes M)$ generated by $\KK(X\otimes M)\cup \{\phi_X(a)\otimes \id_M:a\in A\}$, since $\theta_{x_1\langle m_1,m_2\rangle_A, x_2}\otimes\id_M=\theta_{x_1\otimes m_1,x_2\otimes m_2}$ for all $x_1,x_2\in X$, $m_1,m_2\in M$. 

Besides, if $T\otimes\id_M=0$ for some $T\in \LL(X)$, then $0=\langle Tx\otimes m, Tx\otimes m\rangle = \langle m, \langle Tx,Tx\rangle m\rangle,$  for all $m\in M$, $x\in X$. It follows that $T=0$ because $A$ acts faithfully on $M$.

Notice now that  conjugation by $J$ carries $C$ isomorphically into $L$  because 
\[\begin{array}{c}
J\big(\theta_{x_1\otimes m_1,x_2\otimes m_2}\big)J^{-1}=\theta_{J(x_1\otimes m_1),J(x_2\otimes m_2)},\\
\\
 J(\phi_X(a)\otimes \id_M))J^{-1}=\phi_M(a)\otimes\id_Y,
\end{array}\]
for all $x_1,x_2\in X$, $m_1,m_2\in M$, and $a\in A$. Besides,  $\{\phi_M(a):a\in A\}=\KK(M)$.

Finally, $I(\phi_X(a))=O(\phi_M(a)\otimes \id_Y)=\phi_M(a)\otimes \id_{Y(B)}$, for all $a\in A$.
\end{proof}
\begin{prop}
\label{infbim}
Let $X$, $Y$, and $M$  be  as in Proposition \ref{i1}. Let $\{(X_n,A^X_n,\phi_n^X, \phi_n^A)\}$ and  $\{(M^Y_n,B^Y_n,\phi_n^{M,Y}, \phi_n^B)\}$ be the directed sequences defined in  Example \ref{main}, and let   $(X_\infty, A^X_\infty,\{\lambda^X_n\}, \{\lambda^A_n\})$ and $(M^Y_\infty, B^Y_\infty,\{\mu^M_n\}, \{\mu^B_n\})$, respectively, denote their direct limits.

Then  $M^Y_\infty$ is an $A_\infty^X-B^Y_\infty$ is a Hilbert \cstar-bimodule that  is full on the left. 

Besides, the canonical maps $(\lambda_n^{A}, \mu_n^{M}, \mu_n^B):(A_n^X,M_n^Y, B_n^Y)\rightarrow (A_\infty^X,M_\infty^Y, B_\infty^Y)$ are homomorphisms of Hilbert \cstar-bimodules.

If $M$ is also full on the right, and $Y$ is left non-degenerate as a $B$-module, that is if $\phi_Y(B)Y=Y$, then $M^Y_\infty$ is an $A_\infty^X-B^Y_\infty$ Morita equivalence bimodule.
\end{prop}

\begin{proof}  All the statements involving the right structure  except for the last one, which we discuss at the end,  were taken care of  in Proposition \ref{dlim}, so we focus on the left structure.
We have shown that $M_1^Y$ is an $A_1^X-B_1^Y$ left full Hilbert \cstar-bimodule  by identifying $A_1^X$ with $\KK(M_1^Y)$ via the isomorphism $I$ of Proposition \ref{i1}.

Our aim is to show, in the notation of Example \ref{main},  that $M^Y_n$ is an $A_n^X-B^Y_n$  Hilbert \cstar-bimodule that is full on the left  in a compatible way with the corresponding directed sequences, which  will provide $M^Y_\infty$ with a structure of  $A^X_\infty-B^Y_\infty$ Hilbert \cstar-bimodule. 

First   notice that the map  $J_1:X_1\otimes_{A^X_1}M^Y_1\longrightarrow M^Y_1\otimes_{B^Y_1}Y_1$  given by
\[J_1=(\id_M\otimes L_{B_1^Y,Y}^{-1}\otimes\id_{B^Y_1})(J\otimes\id_{B^Y_1})(\id_X\otimes L_{A_1^X,M^Y_1})\]
is an isomorphism of $A_1^X-B_1^Y$ correspondences. Note that $J_1$ preserves the left action of $A^X_1$ because, by Remark \ref{oops}:
\[\begin{array}{c}
(\id_M\otimes L_{B_1^Y,Y}^{-1}\otimes\id_{B^Y_1})(J\otimes\id_{B^Y_1})(r\otimes \id_{M^Y_1})\\
=[(\id_M\otimes L_{B_1^Y,Y}^{-1})J(r\otimes \id_{M})]\otimes\id_{B^Y_1}\\
= [\big(O(J(r\otimes \id_{M})J^{-1})\otimes \id_Y\big)(\id_M\otimes L_{B_1^Y,Y}^{-1})J]\otimes \id_{B^Y_1}\\
=(I(r)\otimes \id_{Y_1})(\id_M\otimes L_{B_1^Y,Y}^{-1}\otimes\id_{B^Y_1})(J\otimes \id_{B^Y_1}),
\end{array}\]
for $r\in A^X_1$. Besides, $J_1$ is onto and preserves the right action of $B_1^Y$ and the $B_1^Y$-valued  inner product because so do  the maps composed to get $J$.

Notice also that 
\[\xymatrix{
X\otimes_A M\ar[d]_{\phi_0^X\otimes\phi_0^{M,Y}}\ar[r]^J &  M\otimes_B Y\ar[d]^{\phi_0^{M,Y}\otimes\phi_0^Y}\\
X_1\otimes_{A^X_1} M_1^Y\ar[r]^{J_1} &   M_1^Y\otimes_{B^Y_1} Y_1}\]
commutes because 
\[\begin{array}{c}
(\id_M\otimes L_{B_1^Y,Y}^{-1}\otimes\id_{B^Y_1})(J\otimes \id_{B_1^Y})(\id_X\otimes L_{A^X_1, M^Y_1})(\phi_0^X\otimes \phi_0^{M,Y})(xa\otimes mb)\\
=(\id_M\otimes L_{B_1^Y,Y}^{-1}\otimes\id_{B^Y_1})(J\otimes \id_{B_1^Y})(\id_X\otimes L_{A^X_1, M^Y_1})(x\otimes \phi_X(a)\otimes m\otimes \phi_Y(b))\\
=(\id_M\otimes L_{B_1^Y,Y}^{-1}\otimes\id_{B^Y_1})(J\otimes \id_{B_1^Y})[x\otimes I(\phi_X(a))(m\otimes \phi_Y (b))]\\
=(\id_M\otimes L_{B_1^Y,Y}^{-1}\otimes\id_{B^Y_1})(J\otimes \id_{B_1^Y})(x\otimes am\otimes \phi_Y (b))\\
=(\id_M\otimes L_{B_1^Y,Y}^{-1}\otimes\id_{B^Y_1})J(xa\otimes m)\otimes \phi_Y(b)\\
=(\phi_0^{M,Y}\otimes \phi_0^Y)J(xa\otimes mb),
\end{array}\]
for $x\in X$, $a\in A$, $m\in M$ and $b\in B$.

Now  this yields, by Proposition \ref{i1},  an isomorphism $I_2:A_2^X\longrightarrow \KK(M_2^Y)$. Furthermore, the diagram
\[\xymatrix{
A_1^X\ar[d]_{\phi_1^A} \ar[r]^{I_1} & {\KK(M_1^Y)} \ar[d]^{\big(\phi^{M,Y}_1\big)_*}\\
A_2^X \ar[r]^{I_2} & {\KK(M_2^Y)}}\]
commutes, since  by  Proposition \ref{i1} we have, for $r\in A_1^X$: 
\[I_2(\phi_1^{A,X}(r))=\phi_{M_1}(r)\otimes \id_{B^Y_2}=I_1(r)\otimes \id_{B^Y_2}=[\big(\phi^{M,Y}_1\big)_*](I_1(r)),\]
the last equality being due to the fact that
\[\theta_{\phi_1^{M,Y}(m_1b_1),\phi_1^{M,Y}(m_2b_2)}=\theta_{m_1\otimes \phi_1^{B,Y}(b_1),m_2\otimes \phi_1^{B,Y}(b_2)}=\theta_{m_1b_1,m_2b_2}\otimes \id_{B_2^Y},\]
for $m_i\in M_1$, $b_i\in B_1$, and $i=1,2$.

It is clear now that, by iterating this construction,  we get isomorphisms $I_n$ such that the diagram 
\[\xymatrix{
A_n^X\ar[d]_{\phi_n^{A,X}} \ar[r]^{I_n} & {\KK(M_n^Y)} \ar[d]^{\big(\phi^{M,Y}_n\big)_*}\\
A_{n+1}^X \ar[r]^{I_{n+1}} & {\KK(M_{n+1}^Y)}}\]
commutes for all $n\geq 0$.

This shows that $\ainf$ is isomorphic to the direct limit of $\{(\KK(M_n^Y),\big(\phi^{M,Y}_n\big)_* )\}$, which by Proposition \ref{comp} is  $(\KK(M^Y_\infty), (\mu^M_n)_*)$. Therefore  $M_\infty^Y$ is an  $A^X_\infty-B^Y_\infty$ Hilbert \cstar-bimodule that is full on the left,  with left structure defined by
\[\lambda^{A}_n(a_n)\mu^{M}_n(m_n):=\mu^{M}_n(a_nm_n),\ \langle \mu^{M}_n(m_n),\mu^{M}_n(m'_n)\rangle_{A_\infty^X}:=\lambda^{A}_n(\langle m_n,m'_n\rangle_{A_n^X}),\]
for $m_n,m'_n\in M_n^Y$ and $a_n\in A^X_n$, where we write, as we will do from now on,  $a_nm_n$ and $\langle m_n,m'_n\rangle_{A_n^X}$ instead of $[I_n(a_n)](m_n)$ and $I^{-1}_n(\theta_{m_n,m'_n})$, respectively.

Notice that the last equality shows that  $(\lambda_n^{A}, \mu_n^{M})$ is a  homomorphism of left Hilbert \cstar-modules.

If $Y$ is non-degenerate on the left, and $M$ is a Morita equivalence bimodule, then $M^Y_1$ is an $A^X_1-B^Y_1$ Morita equivalence bimodule because $\langle m\otimes r, n\otimes s\rangle=r^*\phi_Y(\langle m,n\rangle_R)s$, for $m,n\in M$ and $r,s\in B^Y_1$.
 Therefore, as one sees by taking an approximate identity for $B^Y_1$,  $\langle M^Y_1,M^Y_1\rangle_R$ contains Im $\phi_Y$ and $\phi_Y(B)\KK(Y)$. But non-degeneracy implies that $\phi_Y(B)\KK(Y)=\KK(Y)$ since, given $x,y\in Y$, then
\[\theta_{x,y}=\theta_{\phi_Y(b)x',y}=\phi_Y(b)\theta_{x',y},\]
for some $x'\in Y$ and $b\in B$. Thus we conclude that $M_1$ is full on the right as well.

It will follow by induction that $M_n$ is full on the right for all $n\geq 0$ once we show that $Y_n$ is always non-degenerate on the left as a  $B^Y_n$-module. In fact:
\[\phi_{Y_n}(B^Y_n)Y_n=\phi_{Y_n}(B^Y_n)(Y_{n-1}\otimes B^Y_n)=B^Y_nY_{n-1}\otimes B^Y_n =Y_{n-1}\otimes B^Y_n,\]
since $B^Y_n\supset \KK(Y_{n-1})$.

Finally, we conclude that in that case $M^Y_\infty$ is full on the right because $\langle M^Y_\infty, M^Y_\infty\rangle$ contains $\mu^B_n\big(\langle M^Y_n, M^Y_n\rangle\big)$ for all $n\geq 0$.
\end{proof}
\begin{rk}
\label{essential}
{\rm{Let  $(Y,\phi_Y)$ be  a correspondence  over a \cstar-algebra $B$, and let $Y_n$, $B^Y_n$, and $M_n^Y$   be as in Example \ref{main}. The proof of Proposition   \ref{infbim} shows that
\begin{enumerate}
\item The $B^Y_n$-left module $Y_n$ is  non-degenerate  for all $n\geq 1$.  Of course, this might fail for $n=0$.
\item If $X$, $Y$ and $M$ are as in Proposition \ref{infbim}, $M$ is full on the right, and $Y$ is non-degenerate, then $M_n^Y$ is an $A_n^X-B^Y_n$ Morita equivalence bimodule such that the $A_n^X-B^Y_n$ correspondences $X_n\otimes M_n^Y$ and   $M_n^Y\otimes Y_n$ are isomorphic for all $n\geq 0$.
\end{enumerate}
}}
\end{rk}

\begin{thm}
\label{meq}
Let $(X,\phi_X)$ and $(Y,\phi_Y)$ be correspondences over the \cstar-algebras $A$ and $B$, respectively. If, in the notation of Example \ref{main}, there exists an $A^X_{n_0}-B^Y_{m_0}$ Morita equivalence bimodule $M$ such that $X_{n_0}\otimes M$ and $M\otimes Y_{m_0}$ are isomorphic  as $A^X_{n_0}-B^Y_{m_0}$ correspondences for some $n_0\geq 0,\  m_0\geq 1$, then the augmented Cuntz-Pimsner \cstar-algebras $\tilde{{\mathcal{O}}}_X$ and $\tilde{{\mathcal{O}}}_Y$ are Morita equivalent.

\end{thm}

\begin{proof}

The bimodules $X_\infty$ and $Y_{\infty}$ and the \cstar-algebras $A^X_\infty$ and $B^Y_\infty$ of Example \ref{main} can be obtained as the limits of the corresponding directed sequences starting, respectively,  at $n_0$ and $m_0$. Besides, the directed sequence $\{(M_n^Y,\phi_n^{M,Y})\}_{n\geq m_0}$ can be constructed as in Example \ref{main}. 
Our aim is to show that $\xinf\otimes_{\ainf}M^Y_\infty$ and $M_\infty^Y\otimes_{B^Y_\infty}Y_\infty$ are isomorphic as $A^X_{\infty}-B^Y_\infty$ Hilbert \cstar-bimodules. It follows from the remarks above that  we can assume that $n_0=m_0=0$ and, in view of the last part of Remark \ref{essential}, that $Y$ is left non-degenerate over $B$.
The result  will then follow from Theorem \ref{eqcp},   Proposition \ref{infbim},   and \cite[4.2]{aee}.

As in Proposition \ref{infbim} and  Example \ref{main} we have the  commuting diagrams:
\[\begin{array}{cc}
\xymatrix{
X_n\ar[r]^{\phi_n^X} \ar[d]_{\lambda^X_n}& X_{n+1}\ar[dl]^{\lambda_{n+1}^X}\\
X_\infty}& \xymatrix{
A^X_n\ar[r]^{\phi_n^{A,X}} \ar[d]_{\lambda^A_n}&A^X_{n+1}\ar[dl]^{\lambda_{n+1}^A}\\
{A^X_\infty}}\\
\xymatrix{
Y_n\ar[r]^{\phi_n^Y} \ar[d]_{\mu^Y_n}&Y_{n+1}\ar[dl]^{\mu_{n+1}^Y}\\
Y_\infty}&\xymatrix{
B^Y_n\ar[r]^{\phi_n^{B,Y}} \ar[d]_{\mu^B_n}&B^Y_{n+1}\ar[dl]^{\mu_{n+1}^B}\\
B^Y_\infty}\\
\xymatrix{
M^Y_n\ar[r]^{\phi_n^{M,Y}} \ar[d]_{\mu^M_n}&M^Y_{n+1}\ar[dl]^{\mu_{n+1}^M}\\
{M^Y_\infty}}
 & \xymatrix{
X_n\otimes_{A^X_n} M^Y_n\ar[d]_{\phi_n^X\otimes\phi_n^{M,Y}}\ar[r]^{J_n} &  M^Y_n\otimes_{B^Y_n} Y_n\ar[d]^{\phi_n^{M,Y}\otimes\phi_n^Y}\\
X_{n+1}\otimes_{A^X_{n+1}} M^Y_{n+1}\ar[r]^{J_{n+1}} &   M^Y_{n+1}\otimes_{B^Y_{n+1}} Y_{n+1}}
\end{array}\]
Notice that, if $m,m'\in M^Y_n$, $y,y'\in Y_n$, then by Propositions \ref{dlim} and \ref{infbim}
\[\begin{array}{ll}
\langle \mu_n^M(m)\otimes \mu_n^Y(y),\mu_n^M(m')\otimes \mu_n^Y(y')\rangle&=\langle  \mu_n^Y(y),\mu_n^B(\langle m,m'\rangle)\mu_n^Y(y')\rangle\\
&=\mu_n^B(\langle y,\langle m,m'\rangle y'\rangle)\\
&=\mu_n^B(\langle m\otimes y, m'\otimes y'\rangle)
\end{array}\]
and 
\[\begin{array}{ll}
\langle \lambda_n^X(x)\otimes \mu_n^M(m),\lambda_n^X(x')\otimes\mu_n^M(m')\rangle&=\langle  \mu_n^M(m),\lambda_n^A(\langle x,x'\rangle)\mu_n^M(m')\rangle\\
&=\langle \mu_n^M(m),\mu_n^M(\langle x,x'\rangle m')\rangle\\
&=\mu_n^B(\langle x\otimes m,x'\otimes m'\rangle)
\end{array}\]
for $x,x'\in X_n$ and $m,m'\in M^Y_n$.

We now want to define $J_\infty:X_\infty\otimes_{A^X_\infty} M^Y_{\infty}\longrightarrow M^Y_{\infty}\otimes_{B^Y_\infty}Y_\infty$ by 
\[J_\infty((\lambda^X_n\otimes \mu_n^M)(x_n\otimes m_n)):=(\mu_n^M\otimes\mu_n^Y)J_n(x_n\otimes m_n).\]

Now, 
\[\langle (\mu_n^M\otimes\mu_n^Y)J_n(x_n\otimes m_n), (\mu_n^M\otimes\mu_n^Y)J_n(x'_n\otimes m'_n)\rangle \]
\[=\mu_n^B(\langle J_n(x_n\otimes m_n),J_n(x'_n\otimes m'_n)\rangle)=\mu_n^B(\langle x_n\otimes m_n,x'_n\otimes m'_n\rangle)\]
\[=\langle (\lambda_n^X\otimes \mu_n^M)(x_n\otimes m_n), (\lambda^X_n\otimes \mu^M_n)(x'_n\otimes  m'_n)\rangle\]
This shows that $J_\infty$ as defined above extends to a right Hilbert \cstar-module homomorphism that  preserves the left  action of $A_\infty$. In fact,  by Proposition \ref{infbim},  given $a_n\in A_n$, $x_n\in X_n$ and $m_n\in M_n$, we have
\[\begin{array}{ll}
J_\infty[\lambda_n(a_n)\cdot(\lambda^X_n(x_n)\otimes \mu_n^M(m_n))]&=J_\infty[\lambda^X_n(a_nx_n)\otimes \mu_n^M(m_n)]\\
&=(\mu^M_n\otimes \mu^Y_n)J_n(a_nx_n\otimes m_n)\\
&=(\mu^M_n\otimes \mu^Y_n)(\phi_{M_n}(a_n)\otimes\id_{Y_n})J_n(x_n\otimes m_n)\\
&=\lambda_n(a_n)\cdot J_\infty(x_n\otimes m_n).
\end{array}\]
Analogous computations show that $J_\infty$ preserves the right action. Besides,  $J_\infty$ is onto because its image contains $\bigcup_n\big(\mu_n^M(M_n)\otimes \mu_n^Y(Y_n)\big)$,  which is dense in $M^Y_\infty\otimes Y_\infty$.

It remains to show that $J_\infty$ preserves the left inner product. This follows as in \cite[1.2]{pic}: if $\xi_0,\xi_1,\xi_2\in X\otimes M$, then
\[\begin{array}{ll}
\langle J_\infty(\xi_0),J_\infty(\xi_1)\rangle_{A_\infty^X}J_\infty(\xi_2)&= J_\infty(\xi_0)\langle J_\infty(\xi_1),J_\infty(\xi_2)\rangle_{B_\infty^Y}\\
&= J_\infty(\xi_0)\langle \xi_1,\xi_2\rangle_{B_\infty^Y}\\
&= J_\infty\big(\xi_0\langle \xi_1,\xi_2\rangle_{B_\infty}\big)\\
&= J_\infty\big(\langle \xi_0,\xi_1\rangle_{A_\infty^X}\xi_2\big)\\
&=\langle \xi_0,\xi_1\rangle_{A_\infty^X}J_\infty(\xi_2)
\end{array}\]
\end{proof}

\begin{rk}
\label{muso}
{\rm{ A similar result was shown  by P. Muhly and B. Solel (\cite{ms}) for  Cuntz-Pimsner \cstar-algebras ${\mathcal {O}}_X$ of correspondences $(X,\phi_X)$ such that $\phi_X$ is injective and $X$ left non-degenerate. Our result for augmented Cuntz-Pimsner \cstar-algebras does not require the faithfulness of $\phi$. Non-degeneracy, however, might play a role, as the following Corollary shows. }}

\end{rk}  

\begin{clly}
Let $(X,\phi_X)$ and $(Y,\phi_Y)$ be correspondences over the \cstar-algebras $A$ and $B$, and let $M$ be a Morita equivalence $A-B$ bimodule such that $X\otimes M$ and $M\otimes Y$ are isomorphic as $A-B$ correspondences. If $Y$ is left non-degenerate, then the augmented Cuntz-Pimsner \cstar-algebras $\ox$ and  $\tilde{{\mathcal{O}}}_Y$ are Morita equivalent.

\end{clly}
\begin{proof}

 By Remark \ref{essential} the conditions in Theorem \ref{meq} are then met.
\end{proof}

\noindent

\end{document}